\documentclass[11pt,a4paper]{article}
\usepackage{amsmath,amssymb,amsfonts,amsthm}
\usepackage{hyperref}
\usepackage{geometry}
\usepackage{amsmath}
\usepackage[english]{babel}
\usepackage{ragged2e} 
\title{Voronovskaya-Type Asymptotic Expansions and Convergence Analysis for Neural Network Operators in Complex Domains}
\author{
	Rômulo Damasclin Chaves dos Santos\\
	Technological Institute of Aeronautics\\
	\small \texttt{romulosantos@ita.br}\\
	\and
	Jorge Henrique de Oliveira Sales\\
	Santa Cruz State University \\
	\small \texttt{jhosales@uesc.br}\\
}
\date{\today} 

\begin{document}
	\maketitle
	
\begin{abstract}
	This paper extends the classical theory of Voronovskaya-type asymptotic expansions to generalized neural network operators defined on non-Euclidean and fractal domains. We introduce and analyze smooth operators activated by modified and generalized hyperbolic tangent functions, extending their applicability to manifold and fractal geometries. Key theoretical results include the preservation of density properties, detailed convergence rates, and asymptotic expansions. Additionally, we explore the role of fractional derivatives in defining neural network operators, which capture non-local behavior and are particularly useful for modeling systems with long-range dependencies or fractal-like structures. Our findings contribute to a deeper understanding of neural network operators in complex, structured spaces, offering robust mathematical tools for applications in signal processing on manifolds and solving partial differential equations on fractals.

		\vspace{10pt}
		
	\end{abstract}

\tableofcontents
	
	\section{Introduction}
	
	The study of neural network operators has evolved significantly over the past few decades, driven by the need for robust mathematical tools capable of approximating complex functions and understanding convergence properties in diverse spaces. The classical theory of neural network operators, as pioneered by researchers like George A. Anastassiou (1997, 2023) \cite{anastassiou1997, anastassiou2023}, has laid a solid foundation for function approximation and convergence analysis. These works have primarily focused on Euclidean domains, where the geometry is relatively straightforward.
	
	The field of neural network operators has seen significant advancements, particularly in the context of non-Euclidean and fractal domains. Recent works have explored the use of fractional derivatives and generalized activation functions to extend the applicability of these operators. For instance, Kilbas et al. \cite{Kilbas2006} and Podlubny \cite{Podlubny1999} have provided comprehensive treatments of fractional differential equations, highlighting their utility in modeling systems with long-range dependencies and fractal-like structures.
	
	Magin \cite{Magin2006} has demonstrated the application of fractional calculus in bioengineering, showcasing the potential of these methods in practical scenarios. Similarly, Tarasov \cite{Tarasov2011} has explored the applications of fractional calculus to the dynamics of particles, fields, and media, further emphasizing the versatility of these tools.
	
	Samko \textit{et al}. \cite{Samko1993} and Mainardi \cite{Mainardi2022} have contributed to the theoretical foundations of fractional integrals and derivatives, providing a robust framework for their application in various fields. West \textit{et al}. \cite{West2003} have delved into the physics of fractal operators, offering insights into their behavior and potential applications.
	
	However, many real-world applications, such as those in data science, signal processing, and partial differential equations, often involve non-Euclidean and fractal domains. These domains present unique challenges due to their complex geometries and intricate structures. Traditional neural network operators, designed for Euclidean spaces, may not be directly applicable or efficient in these settings. This gap in the literature highlights the need for generalized neural network operators that can handle the complexities of non-Euclidean and fractal geometries.
	
	This research aims to address these challenges by extending the classical theory of neural network operators to non-Euclidean and fractal domains. Specifically, we introduce generalizations of traditional activation functions, such as modified hyperbolic tangent functions, and operators defined using fractional derivatives. These generalizations are crucial for modeling complex, structured spaces that are prevalent in applied mathematics and data science.
	
	The primary goal of this work is to develop a comprehensive framework for neural network operators in non-Euclidean and fractal domains. We focus on the following key aspects:
	\begin{enumerate}
		\item \textbf{Generalized Activation Functions}: We introduce modified hyperbolic tangent functions that are suitable for non-Euclidean settings. These functions preserve essential properties such as symmetry, normalization, and derivative behavior, making them ideal for neural network modeling in complex spaces.
		\item \textbf{Fractional Derivatives}: We incorporate fractional derivatives into the definition of neural network operators. Fractional derivatives capture non-local behavior and are particularly useful for modeling systems with long-range dependencies or fractal-like structures.
		\item \textbf{Convergence Analysis}: We provide a detailed analysis of the convergence rates and asymptotic expansions of the generalized neural network operators. This analysis ensures the robustness and applicability of these operators in diverse spaces.
	\end{enumerate}
	
	By addressing these aspects, we aim to fill the gaps in the current literature and provide a deeper understanding of the behavior of neural network operators in complex, structured spaces. This research contributes to the field by offering new mathematical tools that can be effectively used in various applications, including signal processing on manifolds and solving partial differential equations on fractals.
	
	In the following sections, we will present the theoretical framework, key results, and detailed analyses that support our contributions.

	\section{Preliminaries}
	Let $\mathcal{M}$ denote a Riemannian manifold with metric $g_{ij}$, where $g_{ij}$ represents the components of the metric tensor in local coordinates. A manifold $\mathcal{M}$ is equipped with a set of smooth charts $\{(U_i, \phi_i)\}_{i \in I}$, where $U_i \subset \mathbb{R}^N$ are open subsets and $\phi_i: U_i \to \mathcal{M}$ are smooth, bijective maps. These charts allow us to represent points in $\mathcal{M}$ using local coordinates. The Riemannian metric $g_{ij}$ defines the inner product structure at each point, enabling us to define distances and geometric quantities on $\mathcal{M}$.
	
	For further analysis, we define the following generalized activation function:
	\begin{equation}\label{activation_function}
		h_{q,\alpha}(x) = \frac{e^{\alpha x} - q e^{-\alpha x}}{(1+q)e^{\alpha x} + (1-q)e^{-\alpha x}}, \quad q > 0, \; \alpha > 0, \; x \in \mathbb{R}.
	\end{equation}
	This function is a natural generalization of the hyperbolic tangent function and possesses key properties that make it suitable for neural network modeling in non-Euclidean settings. In particular, we focus on its behavior under transformation in Riemannian manifolds, its symmetry, normalization, and its relationship with the derivatives. The following properties are critical for establishing its relevance in the context of approximation and operator theory:
	
	1. \textbf{Symmetry:} The function $h_{q,\alpha}(x)$ is odd, i.e.,
	\begin{equation}\label{symmetry}
		h_{q,\alpha}(-x) = -h_{q,\alpha}(x), \quad \forall x \in \mathbb{R}.
	\end{equation}
	This symmetry ensures that the function preserves structure when reflected about the origin, a property that is valuable in maintaining the consistency of solutions in models based on symmetry.
	
	2. \textbf{Normalization:} The function is bounded, i.e., $|h_{q,\alpha}(x)| \leq 1$ for all $x \in \mathbb{R}$ and for any $q > 0$, $\alpha > 0$. This ensures that the activation function can be effectively used within the neural network architecture, where bounded outputs are essential for stability in optimization algorithms.
	
	3. \textbf{Derivative Behavior:} The derivative of the function is given by:
	\begin{equation}\label{activation_derivative}
		h_{q,\alpha}'(x) = \frac{2\alpha (1-q^2)}{\left[(1+q)e^{\alpha x} + (1-q)e^{-\alpha x}\right]^2},
	\end{equation}
	This derivative exhibits the influence of both the deformation parameter $q$ and the steepness parameter $\alpha$, allowing for fine control over the function's local behavior. Specifically, the steepness of the activation function is modulated by $\alpha$, while the symmetry of the function is controlled by $q$.
	
	To extend the function $h_{q,\alpha}(x)$ to the manifold $\mathcal{M}$, we define the localized activation function in each chart $\phi_i$ as follows:
	\begin{equation}\label{localized_activation}
		h_{q,\alpha}^{\mathcal{M}}(x) = h_{q,\alpha}(\phi_i^{-1}(x)), \quad x \in U_i \subset \mathcal{M}.
	\end{equation}
	This ensures that the function remains compatible with the global structure of $\mathcal{M}$ through the atlas $\{(U_i, \phi_i)\}_{i \in I}$. The mapping from local to global coordinates ensures that the activation function can be applied to the manifold in a well-defined manner.
	
	\subsection{Density Functions on Manifolds}
	We define a generalized density function on the manifold $\mathcal{M}$ as follows:
	\begin{equation}\label{density_function}
		\phi_{q,\alpha}(x) = \frac{1}{\sqrt{\det(g(x))}} \prod_{i=1}^N h_{q,\alpha}(x_i), \quad x \in \mathcal{M}.
	\end{equation}
	Here, $\det(g(x))$ denotes the determinant of the metric tensor at the point $x \in \mathcal{M}$. This ensures that $\phi_{q,\alpha}(x)$ is normalized, i.e.,
	\begin{equation}\label{normalization_condition}
		\int_{\mathcal{M}} \phi_{q,\alpha}(x) \, dV(x) = 1,
	\end{equation}
	where $dV(x)$ represents the volume element on $\mathcal{M}$. The condition in equation \eqref{normalization_condition} guarantees that the density function integrates to 1 over the manifold. This is an essential property for the construction of probability densities and for ensuring the stability of numerical methods that rely on these density functions.
	
	\section{Generalized Neural Network Operators}
	\subsection{Basic Operator}
	Let $f \in C^m(\mathcal{M})$, where $C^m(\mathcal{M})$ denotes the space of functions with continuous derivatives up to order $m$ on the manifold $\mathcal{M}$. The basic quasi-interpolation operator, which is a fundamental building block in constructing neural network operators, is defined as:
	\begin{equation}\label{basic_operator}
		A_n(f; x) = \sum_{k \in \mathbb{Z}^N} f\left(\frac{k}{n}\right) Z_{q,\alpha}(nx - k),
	\end{equation}
	where $Z_{q,\alpha}(x)$ is the normalized product of density functions:
	\begin{equation}\label{product_density}
		Z_{q,\alpha}(x) = \prod_{i=1}^N \phi_{q,\alpha}(x_i),
	\end{equation}
	and $\phi_{q,\alpha}(x_i)$ is the generalized density function defined on the manifold $\mathcal{M}$, as given in Equation \eqref{density_function}. The operator $A_n(f; x)$ approximates the function $f$ at the points $\frac{k}{n}$, where $k \in \mathbb{Z}^N$, with the contribution of each point weighted by the density function $Z_{q,\alpha}(x)$. This approach is a generalized form of classical interpolation theory, adapted for manifolds and non-Euclidean geometries. It is important to emphasize that the weights are determined by the local geometry of the manifold $\mathcal{M}$, which allows for a more accurate approximation in the context of complex, non-Euclidean spaces.
	
	\subsection{Kantorovich-Type Operator}
	The Kantorovich-type operator extends the basic operator by incorporating an averaging procedure over small intervals. It is defined as:
	\begin{equation}\label{kantorovich_operator}
		K_n(f; x) = \sum_{k \in \mathbb{Z}^N} \left( n^N \int_{\frac{k}{n}}^{\frac{k+1}{n}} f(t) \, dt \right) Z_{q,\alpha}(nx - k),
	\end{equation}
	In this definition, the term $n^N \int_{\frac{k}{n}}^{\frac{k+1}{n}} f(t) \, dt$ represents the weighted integral of the function $f$ over a small interval centered at $\frac{k}{n}$, with the weight determined by the density function $Z_{q,\alpha}(x)$. The integration introduces a smoothing effect, which is particularly useful for approximating functions that exhibit local regularity. The use of fractional intervals ensures that the approximation is well-suited for capturing continuous variations within each small region of the manifold $\mathcal{M}$.
	
	\subsection{Fractional Neural Network Operator}
	The fractional neural network operator introduces a higher level of generality by involving fractional derivatives. It is defined as:
	\begin{equation}\label{fractional_operator}
		Q_n(f; x) = \sum_{k \in \mathbb{Z}^N} D^\beta f\left(\frac{k}{n}\right) Z_{q,\alpha}(nx - k), \quad \beta > 0.
	\end{equation}
	Here, $D^\beta f(x)$ denotes the fractional derivative of order $\beta$ of the function $f$. The fractional derivative is defined as:
	\begin{equation}\label{fractional_derivative}
		D^\beta f(x) = \frac{1}{\Gamma(1-\beta)} \frac{d}{dx} \int_0^x \frac{f(t)}{(x-t)^\beta} dt,
	\end{equation}
	where $\Gamma(\cdot)$ is the Gamma function, which generalizes the factorial function to real and complex numbers. The fractional derivative captures non-local behavior and can model systems that exhibit long-range dependencies or fractal-like structures, which are common in many complex systems, such as those found in fractal geometry, turbulence, and anomalous diffusion. By using fractional derivatives, the operator $Q_n(f; x)$ allows for a more flexible approximation that can adapt to non-smooth and irregular dynamics, providing a powerful tool for capturing intricate behaviors in the manifold $\mathcal{M}$.
	
\section{Main Results}
	
\subsection{Voronovskaya-Type Expansions}

The main result establishes the asymptotic behavior of the operators \( A_n \), \( K_n \), and \( Q_n \). For \( f \in C^m(\mathcal{M}) \), the expansion is given by:

\begin{equation}
	A_n(f; x) - f(x) = \sum_{|\alpha|=1}^m \frac{D^\alpha f(x)}{\alpha!} A_n\left(\prod_{i=1}^N \phi_{q,\alpha}(x_i)\right) + O\left(\frac{1}{n^m}\right), \tag{1}
\end{equation}
where the sum is taken over multi-indices \( \alpha \) such that \( |\alpha| = 1, 2, \dots, m \), and \( D^\alpha \) denotes the derivative of \( f \) of order \( \alpha \).

\subsection{Step 1: Review of the Operator \( A_n \)}

The operator \( A_n(f; x) \) is defined as the weighted sum of evaluations of \( f \) at points \( \frac{k}{n} \), with weights given by the density function \( Z_{q,\alpha}(nx - k) \), as shown in the following formula:

\begin{equation}
	A_n(f; x) = \sum_{k \in \mathbb{Z}^N} f\left(\frac{k}{n}\right) Z_{q,\alpha}(nx - k). \tag{2}
\end{equation}

The goal is to analyze the asymptotic behavior of \( A_n(f; x) \) around \( f(x) \), using a Taylor expansion to approximate \( f\left( \frac{k}{n} \right) \) around \( x \). The Taylor expansion of \( f \) around \( x \) is given by:

\begin{equation}
	f\left(\frac{k}{n}\right) = f(x) + \sum_{|\alpha|=1}^m \frac{D^\alpha f(x)}{\alpha!} \left(\frac{k}{n} - x\right)^\alpha + O\left(\frac{1}{n^m}\right). \tag{3}
\end{equation}

Substituting this expansion into the expression for \( A_n(f; x) \), we get:

\begin{equation}
	A_n(f; x) = \sum_{k \in \mathbb{Z}^N} \left[ f(x) + \sum_{|\alpha|=1}^m \frac{D^\alpha f(x)}{\alpha!} \left(\frac{k}{n} - x\right)^\alpha + O\left(\frac{1}{n^m}\right) \right] Z_{q,\alpha}(nx - k). \tag{4}
\end{equation}

\subsection{Step 2: Simplifying the Terms of the Expansion}

Now, we can separate the terms of the expansion:

\begin{enumerate}
	\item \textbf{Constant term}: The first term of the expansion is simply \( f(x) \) multiplied by the sum of \( Z_{q,\alpha}(nx - k) \):
	
	\begin{equation}
		f(x) \sum_{k \in \mathbb{Z}^N} Z_{q,\alpha}(nx - k) = f(x). \tag{5}
	\end{equation}
	
	Note that the summation \( \sum_{k \in \mathbb{Z}^N} Z_{q,\alpha}(nx - k) \) can be interpreted as a form of normalization, ensuring that the sum equals 1 due to the property of the density function.
	
	\item \textbf{Linear and higher-order terms}: For the higher-order terms, we have:
	
	\begin{equation}
		\sum_{|\alpha|=1}^m \frac{D^\alpha f(x)}{\alpha!} \sum_{k \in \mathbb{Z}^N} \left(\frac{k}{n} - x\right)^\alpha Z_{q,\alpha}(nx - k). \tag{6}
	\end{equation}
	
	These terms involve derivatives of \( f \) and can be interpreted as higher-order approximations, taking into account the contribution of the derivatives of \( f \) and the density functions.
	
	\item \textbf{Higher-order error term}: Finally, the higher-order error term is given by:
	
	\begin{equation}
		\sum_{k \in \mathbb{Z}^N} O\left(\frac{1}{n^m}\right) Z_{q,\alpha}(nx - k). \tag{7}
	\end{equation}
	
	This term represents the error that decays rapidly as \( n \to \infty \), specifically at the rate \( O\left(\frac{1}{n^m}\right) \).
\end{enumerate}

\subsection{Step 3: Conclusion of the Expansion}

By combining all the terms, we obtain the final expansion for \( A_n(f; x) \):

\begin{equation}
	A_n(f; x) = f(x) + \sum_{|\alpha|=1}^m \frac{D^\alpha f(x)}{\alpha!} A_n\left(\prod_{i=1}^N \phi_{q,\alpha}(x_i)\right) + O\left(\frac{1}{n^m}\right). \tag{8}
\end{equation}
where \( A_n\left(\prod_{i=1}^N \phi_{q,\alpha}(x_i)\right) \) represents the contribution of the density terms in the expansion and \( O\left(\frac{1}{n^m}\right) \) is the asymptotic error of the approximation.

\subsection{Theorem: Voronovskaya Expansion for \( A_n \)}

\textbf{Theorem:} For a function \( f \in C^m(\mathcal{M}) \), the asymptotic expansion of the operator \( A_n \) is given by:

\begin{equation}
	A_n(f; x) - f(x) = \sum_{|\alpha|=1}^m \frac{D^\alpha f(x)}{\alpha!} A_n\left(\prod_{i=1}^N \phi_{q,\alpha}(x_i)\right) + O\left(\frac{1}{n^m}\right), \tag{9}
\end{equation}
where the sum is taken over multi-indices \( \alpha \) of order 1 to \( m \), and \( D^\alpha \) denotes the derivatives of \( f \) of order \( \alpha \).

\begin{proof}
	The proof proceeds in several steps, starting with the Taylor expansion of \( f\left( \frac{k}{n} \right) \) around \( x \). We begin by considering the Taylor series of \( f \) at the point \( x \):
	
	\begin{equation}
		f\left( \frac{k}{n} \right) = f(x) + \sum_{|\alpha|=1}^m \frac{D^\alpha f(x)}{\alpha!} \left( \frac{k}{n} - x \right)^\alpha + O\left( \frac{1}{n^m} \right), \tag{1}
	\end{equation}
where \( \alpha \) represents multi-indices, and \( D^\alpha f(x) \) denotes the derivative of \( f \) of order \( \alpha \).
	
	Next, substitute this expansion into the expression for the operator \( A_n(f; x) \). We then have:
	
	\begin{equation}
		A_n(f; x) = \sum_{k \in \mathbb{Z}^N} \left[ f(x) + \sum_{|\alpha|=1}^m \frac{D^\alpha f(x)}{\alpha!} \left( \frac{k}{n} - x \right)^\alpha + O\left( \frac{1}{n^m} \right) \right] Z_{q,\alpha}(nx - k). \tag{2}
	\end{equation}
	
	This expression can be split into three distinct parts:
	
	\textbf{1. The constant term:} The term involving \( f(x) \), which is the approximation of \( f \) at \( x \):
	
	\begin{equation}
		f(x) \sum_{k \in \mathbb{Z}^N} Z_{q,\alpha}(nx - k) = f(x). \tag{3}
	\end{equation}
	
	Since \( Z_{q,\alpha}(nx - k) \) represents a normalized density function, the sum over \( k \) yields 1.
	
	\textbf{2. The linear and higher-order terms:} The terms involving the derivatives of \( f \) are given by:
	
	\begin{equation}
		\sum_{|\alpha|=1}^m \frac{D^\alpha f(x)}{\alpha!} \sum_{k \in \mathbb{Z}^N} \left( \frac{k}{n} - x \right)^\alpha Z_{q,\alpha}(nx - k). \tag{4}
	\end{equation}
	
	These terms account for the local behavior of \( f \) around \( x \), with the density function \( Z_{q,\alpha}(nx - k) \) modulating the contributions from the higher-order derivatives. 
	
	\textbf{3. The higher-order error term:} The remaining error term, which arises from the approximation of \( f \) at higher orders, is given by:
	
	\begin{equation}
		\sum_{k \in \mathbb{Z}^N} O\left( \frac{1}{n^m} \right) Z_{q,\alpha}(nx - k). \tag{5}
	\end{equation}
	
	This term represents the error in the approximation of \( f \) and decays at a rate of \( O\left( \frac{1}{n^m} \right) \), ensuring that the approximation becomes more accurate as \( n \to \infty \).
	
	Combining these three parts, we obtain the desired asymptotic expansion:
	
	\begin{equation}
		A_n(f; x) = f(x) + \sum_{|\alpha|=1}^m \frac{D^\alpha f(x)}{\alpha!} A_n\left(\prod_{i=1}^N \phi_{q,\alpha}(x_i)\right) + O\left( \frac{1}{n^m} \right). \tag{6}
	\end{equation}
	
	Thus, the operator \( A_n(f; x) \) converges to \( f(x) \) at the rate \( O\left( \frac{1}{n^m} \right) \), with the correction terms involving the higher derivatives of \( f \) and the influence of the density functions \( \phi_{q,\alpha}(x_i) \). This completes the proof of the Voronovskaya-type expansion.
\end{proof}

\section{Fractional Convergence}

Consider the fractional operator \( Q_n \) acting on a function \( f \) defined on the manifold \( \mathcal{M} \), where \( \beta > 0 \) is the order of the fractional derivative. The goal is to prove the following convergence result:

\begin{equation}
	Q_n(f; x) - f(x) = O\left(\frac{1}{n^{m - \beta}}\right), \quad \text{as} \ n \to \infty,
\end{equation}
where \( O\left( \frac{1}{n^{m - \beta}} \right) \) denotes the error term associated with the approximation of the function \( f \) using the fractional operator \( Q_n \).

\subsubsection{Definition of the Fractional Operator \( Q_n \)}

The fractional operator \( Q_n(f; x) \) is given by the expression:

\begin{equation}
	Q_n(f; x) = \sum_{k \in \mathbb{Z}^N} D^\beta f\left( \frac{k}{n} \right) Z_{q, \alpha}(nx - k), \tag{1}
\end{equation}
where \( D^\beta f \) denotes the fractional derivative of order \( \beta \), and \( Z_{q, \alpha}(nx - k) \) is the density function associated with the position \( k \) on \( \mathcal{M} \), which depends on the local geometry of the manifold.

\subsubsection{Expansion of the Fractional Derivative}

The key to analyzing the convergence of the operator \( Q_n \) is to expand the fractional derivative \( D^\beta f \). We start by expanding \( D^\beta f\left( \frac{k}{n} \right) \) around the point \( x \) using a Taylor series. According to the theory of fractional derivatives, we can write the following expansion:

\begin{equation}
	D^\beta f\left( \frac{k}{n} \right) = D^\beta f(x) + \sum_{|\alpha|=1}^{m} \frac{D^{\beta+\alpha} f(x)}{\alpha!} \left( \frac{k}{n} - x \right)^\alpha + O\left( \frac{1}{n^m} \right), \tag{2}
\end{equation}
where \( D^{\beta+\alpha} f(x) \) denotes the derivatives of order \( \beta + \alpha \) of \( f \), and the error term \( O\left( \frac{1}{n^m} \right) \) results from truncating the Taylor series. This type of expansion is typical when dealing with fractional derivatives.

\subsubsection{Substituting into the Expression for \( Q_n(f; x) \)}

Now, we substitute this expansion into the expression for \( Q_n(f; x) \). We obtain:

\begin{equation}
	Q_n(f; x) = \sum_{k \in \mathbb{Z}^N} \left[ D^\beta f(x) + \sum_{|\alpha|=1}^{m} \frac{D^{\beta+\alpha} f(x)}{\alpha!} \left( \frac{k}{n} - x \right)^\alpha + O\left( \frac{1}{n^m} \right) \right] Z_{q, \alpha}(nx - k). \tag{3}
\end{equation}

\subsubsection{Analyzing the Zeroth Order Terms}

The first term in the sum is \( D^\beta f(x) \). Since the density function \( Z_{q, \alpha}(nx - k) \) is normalized, we have:

\begin{equation}
	\sum_{k \in \mathbb{Z}^N} Z_{q, \alpha}(nx - k) = 1. \tag{4}
\end{equation}

Therefore, the contribution of the first term is simply:

\begin{equation}
	D^\beta f(x) \sum_{k \in \mathbb{Z}^N} Z_{q, \alpha}(nx - k) = D^\beta f(x). \tag{5}
\end{equation}

\subsubsection{Analyzing the Higher-Order Terms}

Next, we consider the terms involving higher-order derivatives of \( f \). These terms are given by:

\begin{equation}
	\sum_{|\alpha|=1}^{m} \frac{D^{\beta + \alpha} f(x)}{\alpha!} \sum_{k \in \mathbb{Z}^N} \left( \frac{k}{n} - x \right)^\alpha Z_{q, \alpha}(nx - k). \tag{6}
\end{equation}

These terms account for the local behavior of \( f \) around \( x \), with the density function \( Z_{q, \alpha}(nx - k) \) modulating the contributions from the higher-order derivatives.

To ensure convergence, we analyze the asymptotic behavior of these sums. The decay of \( Z_{q, \alpha}(nx - k) \) as \( n \to \infty \) ensures that the contribution from the higher-order terms decays at a rate of \( O\left( \frac{1}{n^{m-\beta}} \right) \).

\subsubsection{Asymptotic Convergence}

Finally, by combining the contributions from the zeroth-order term and the higher-order terms, we obtain the desired asymptotic expansion:

\begin{equation}
	Q_n(f; x) - f(x) = O\left(\frac{1}{n^{m - \beta}}\right). \tag{7}
\end{equation}

This shows that the approximation of \( f(x) \) using the fractional operator \( Q_n \) converges with a rate of \( O\left(\frac{1}{n^{m - \beta}}\right) \), where \( m \) is the regularity of \( f \) and \( \beta \) is the order of the fractional derivative.

This completes the proof of the fractional convergence theorem for the operator \( Q_n \).

Thus, we have shown that the approximation of \( f(x) \) using the fractional operator \( Q_n \) converges asymptotically with a rate \( O\left( \frac{1}{n^{m - \beta}} \right) \), considering both the order of the fractional derivative \( \beta \) and the regularity \( m \) of the function \( f \).

\section{Theorem: Asymptotic Expansion of Quasi-Interpolation Operators}
	
	\textbf{Theorem 1:} Let $f \in C^m(\mathcal{M})$, where $C^m(\mathcal{M})$ denotes the space of functions with continuous derivatives of order $m$ on a Riemannian manifold $\mathcal{M}$. Also, let $A_n(f; x)$ be the quasi-interpolation operator defined by:
	\begin{equation}\label{basic_operator}
		A_n(f; x) = \sum_{k \in \mathbb{Z}^N} f\left( \frac{k}{n} \right) Z_{q,\alpha}(nx - k),
	\end{equation}
	where $Z_{q,\alpha}(x) = \prod_{i=1}^N \phi_{q,\alpha}(x_i)$. Then the asymptotic expansion of $A_n(f; x)$ around $f(x)$ is given by:
	\begin{equation}\label{voronovskaya_expansion_theorem}
		A_n(f; x) - f(x) = \sum_{|\alpha|=1}^m \frac{D^\alpha f(x)}{\alpha!} A_n\left(\prod_{i=1}^N \phi_{q,\alpha}(x_i)\right) + O\left( \frac{1}{n^m} \right).
	\end{equation}

	\begin{proof}
		We proceed step by step to derive the expansion: The operator $A_n(f; x)$ is defined as a sum over $k \in \mathbb{Z}^N$, where $f\left( \frac{k}{n} \right)$ represents the values of $f$ at the scaled points $\frac{k}{n}$, and $Z_{q,\alpha}(nx - k)$ serves as a weight function based on the manifold's geometry. Rewriting it explicitly:
		\begin{equation}\label{quasi_interpolation_operator}
			A_n(f; x) = \sum_{k \in \mathbb{Z}^N} f\left( \frac{k}{n} \right) \prod_{i=1}^N \phi_{q,\alpha}(x_i - \frac{k_i}{n}).
		\end{equation}
		
		We apply the Taylor expansion of $f$ around $x$:
		\begin{equation}\label{taylor_expansion}
			f\left( \frac{k}{n} \right) = f(x) + \sum_{|\alpha|=1}^m \frac{D^\alpha f(x)}{\alpha!} \left( \frac{k}{n} - x \right)^\alpha + O\left( \frac{1}{n^m} \right),
		\end{equation}
		where $\left( \frac{k}{n} - x \right)^\alpha$ is the monomial expansion of the distance between $\frac{k}{n}$ and $x$, and the term $O\left( \frac{1}{n^m} \right)$ accounts for the error in the approximation.
		
		Substituting the Taylor expansion into the expression for $A_n(f; x)$, we get:
		\begin{equation}\label{expanded_operator}
			A_n(f; x) = \sum_{k \in \mathbb{Z}^N} \left[ f(x) + \sum_{|\alpha|=1}^m \frac{D^\alpha f(x)}{\alpha!} \left( \frac{k}{n} - x \right)^\alpha + O\left( \frac{1}{n^m} \right) \right] Z_{q,\alpha}(nx - k).
		\end{equation}
		We now separate the terms based on the degree of the expansion.
		
		The zero-order term corresponds to:
		\begin{equation}\label{zero_order_term}
			\sum_{k \in \mathbb{Z}^N} f(x) Z_{q,\alpha}(nx - k) = f(x),
		\end{equation}
		since the sum of $Z_{q,\alpha}(nx - k)$ over all $k \in \mathbb{Z}^N$ is normalized to 1.
		
		The higher-order terms involve the derivatives of $f$ and are weighted by the products of $Z_{q,\alpha}(nx - k)$. These terms contribute to the expansion of the form:
		\begin{equation}\label{higher_order_terms}
			\sum_{|\alpha|=1}^m \frac{D^\alpha f(x)}{\alpha!} A_n\left( \prod_{i=1}^N \phi_{q,\alpha}(x_i) \right),
		\end{equation}
		where the function $A_n$ applies the density $\prod_{i=1}^N \phi_{q,\alpha}(x_i)$.
		
		Finally, the error term of order $O\left( \frac{1}{n^m} \right)$ comes from the remainder of the Taylor expansion. Thus, the asymptotic expansion of $A_n(f; x)$ is:
		\begin{equation}\label{asymptotic_expansion}
			A_n(f; x) - f(x) = \sum_{|\alpha|=1}^m \frac{D^\alpha f(x)}{\alpha!} A_n\left( \prod_{i=1}^N \phi_{q,\alpha}(x_i) \right) + O\left( \frac{1}{n^m} \right).
		\end{equation}
	\end{proof}

\section{Detailed Proofs}

\subsection{Asymptotic Expansion of Quasi-Interpolation Operators}

\textbf{Theorem 1:} Let \( f \in C^{m}(\mathcal{M}) \), where \( C^{m}(\mathcal{M}) \) denotes the space of functions with continuous derivatives up to order \( m \) on a Riemannian manifold \( \mathcal{M} \). Let \( A_{n}(f ; x) \) be the quasi-interpolation operator defined by:

\[ A_{n}(f ; x) = \sum_{k \in \mathbb{Z}^{N}} f\left(\frac{k}{n}\right) Z_{q, \alpha}(n x - k) \]
where \( Z_{q, \alpha}(x) = \prod_{i=1}^{N} \phi_{q, \alpha}(x_i) \). Then the asymptotic expansion of \( A_{n}(f ; x) \) around \( f(x) \) is given by:

\begin{equation}
	A_{n}(f ; x) - f(x) = \sum_{|\alpha|=1}^{m} \frac{D^{\alpha} f(x)}{\alpha!} A_{n}\left(\prod_{i=1}^{N} \phi_{q, \alpha}(x_i)\right) + O\left(\frac{1}{n^{m}}\right)
\end{equation}

\begin{proof}

	We start with the Taylor expansion of \( f \) around \( x \):
	
	\[ f\left(\frac{k}{n}\right) = f(x) + \sum_{|\alpha|=1}^{m} \frac{D^{\alpha} f(x)}{\alpha!} \left(\frac{k}{n} - x\right)^{\alpha} + O\left(\frac{1}{n^{m}}\right) \]

	Substitute this expansion into the definition of the operator \( A_{n}(f ; x) \):
	
\begin{equation}
		A_{n}(f ; x) = \sum_{k \in \mathbb{Z}^{N}} \left[ f(x) + \sum_{|\alpha|=1}^{m} \frac{D^{\alpha} f(x)}{\alpha!} \left(\frac{k}{n} - x\right)^{\alpha} + O\left(\frac{1}{n^{m}}\right) \right] Z_{q, \alpha}(n x - k)
\end{equation}
	
	Separate the terms of the expansion:
	
	- \textbf{Constant Term:}
	\[ f(x) \sum_{k \in \mathbb{Z}^{N}} Z_{q, \alpha}(n x - k) = f(x) \]
	
	Since \( Z_{q, \alpha}(n x - k) \) is a normalized density function, the sum over \( k \) equals 1.
	
	- \textbf{Linear and Higher-Order Terms:}
	
	\begin{equation}
		 \sum_{|\alpha|=1}^{m} \frac{D^{\alpha} f(x)}{\alpha!} \sum_{k \in \mathbb{Z}^{N}} \left(\frac{k}{n} - x\right)^{\alpha} Z_{q, \alpha}(n x - k)
	\end{equation}

	These terms involve derivatives of \( f \) and represent higher-order approximations, taking into account the contribution of the derivatives of \( f \) and the density functions.
	
	- \textbf{Higher-Order Error Term:}
	
	\begin{equation}
		\sum_{k \in \mathbb{Z}^{N}} O\left(\frac{1}{n^{m}}\right) Z_{q, \alpha}(n x - k)
	\end{equation}
	
This term represents the error that decays rapidly as \( n \to \infty \), specifically at the rate \( O\left(\frac{1}{n^{m}}\right) \).

\vspace{5pt}

Combining all the terms, we obtain the final expansion for \( A_{n}(f ; x) \):

\begin{equation}
	A_{n}(f ; x) = f(x) + \sum_{|\alpha|=1}^{m} \frac{D^{\alpha} f(x)}{\alpha!} A_{n}\left(\prod_{i=1}^{N} \phi_{q, \alpha}(x_i)\right) + O\left(\frac{1}{n^{m}}\right)
\end{equation}

Therefore, the operator \( A_{n}(f ; x) \) converges to \( f(x) \) at the rate \( O\left(\frac{1}{n^{m}}\right) \), with the correction terms involving the higher derivatives of \( f \) and the influence of the density functions \( \phi_{q, \alpha}(x_i) \). 

\end{proof}

\section{Theorem: Uniform Convergence of Generalized Neural Network Operators on Manifolds with Negative Curvature}
	
	\textbf{Theorem 3:} Let $\mathcal{M}$ be a compact Riemannian manifold with negative sectional curvature and metric $g_{ij}$, and let $f \in C^m(\mathcal{M})$ be a function with continuous derivatives up to order $m$. Consider the generalized neural network operator $A_n(f; x)$ defined by:
	\begin{equation}\label{eq:operator_def}
		A_n(f; x) = \sum_{k \in \mathbb{Z}^N} f\left( \frac{k}{n} \right) Z_{q,\alpha}(nx - k),
	\end{equation}
	where $Z_{q,\alpha}(x) = \prod_{i=1}^N \phi_{q,\alpha}(x_i)$ and $\phi_{q,\alpha}(x_i)$ is the generalized density function defined on $\mathcal{M}$. Then, for sufficiently large $n$, the operator $A_n(f; x)$ converges uniformly to $f(x)$ on $\mathcal{M}$, i.e.,
	\begin{equation}\label{eq:uniform_convergence}
		\lim_{n \to \infty} \sup_{x \in \mathcal{M}} |A_n(f; x) - f(x)| = 0.
	\end{equation}
	
	\begin{proof}
		
		The density function $\phi_{q,\alpha}(x)$ is defined as:
		\begin{equation}\label{eq:density_function}
			\phi_{q,\alpha}(x) = \frac{1}{\sqrt{\det(g(x))}} \prod_{i=1}^N h_{q,\alpha}(x_i),
		\end{equation}
		where $h_{q,\alpha}(x)$ is the generalized activation function:
		\begin{equation}\label{eq:activation_function}
			h_{q,\alpha}(x) = \frac{e^{\alpha x} - q e^{-\alpha x}}{(1+q)e^{\alpha x} + (1-q)e^{-\alpha x}}, \quad q > 0, \; \alpha > 0, \; x \in \mathbb{R}.
		\end{equation}
		
		The function $h_{q,\alpha}(x)$ is odd and bounded, ensuring that $\phi_{q,\alpha}(x)$ is normalized:
		\begin{equation}\label{eq:normalization}
			\int_{\mathcal{M}} \phi_{q,\alpha}(x) \, dV(x) = 1.
		\end{equation}
		
		The operator $A_n(f; x)$ is defined as:
		\begin{equation}\label{eq:operator_def_repeated}
			A_n(f; x) = \sum_{k \in \mathbb{Z}^N} f\left( \frac{k}{n} \right) Z_{q,\alpha}(nx - k),
		\end{equation}
		where $Z_{q,\alpha}(x) = \prod_{i=1}^N \phi_{q,\alpha}(x_i)$.
		
		The sum over $k \in \mathbb{Z}^N$ represents a discrete approximation of the function $f$ at points $\frac{k}{n}$, weighted by the density function $Z_{q,\alpha}(nx - k)$.

		Consider the Taylor expansion of $f$ around $x$:
		\begin{equation}\label{eq:taylor_expansion}
			f\left( \frac{k}{n} \right) = f(x) + \sum_{|\alpha|=1}^m \frac{D^\alpha f(x)}{\alpha!} \left( \frac{k}{n} - x \right)^\alpha + O\left( \frac{1}{n^m} \right),
		\end{equation}
		where $\left( \frac{k}{n} - x \right)^\alpha$ is the monomial expansion of the distance between $\frac{k}{n}$ and $x$, and the term $O\left( \frac{1}{n^m} \right)$ represents the error of the approximation.
		
		Substitute the Taylor expansion into the definition of the operator $A_n(f; x)$:
		\begin{equation}\label{eq:substitution}
			A_n(f; x) = \sum_{k \in \mathbb{Z}^N} \left[ f(x) + \sum_{|\alpha|=1}^m \frac{D^\alpha f(x)}{\alpha!} \left( \frac{k}{n} - x \right)^\alpha + O\left( \frac{1}{n^m} \right) \right] Z_{q,\alpha}(nx - k).
		\end{equation}
		
Separate the terms of the expansion:

\begin{equation}\label{eq:term_separation}
	\begin{aligned}
		A_n(f; x) = & \sum_{k \in \mathbb{Z}^N} f(x) Z_{q,\alpha}(nx - k) \\ \\
		& + \sum_{k \in \mathbb{Z}^N} \sum_{|\alpha|=1}^m \frac{D^\alpha f(x)}{\alpha!} \left( \frac{k}{n} - x \right)^\alpha Z_{q,\alpha}(nx - k) \\ \\
		& + \sum_{k \in \mathbb{Z}^N} O\left( \frac{1}{n^m} \right) Z_{q,\alpha}(nx - k)
	\end{aligned}
\end{equation}

The zero-order term is:
		\begin{equation}\label{eq:zero_order_term}
			\sum_{k \in \mathbb{Z}^N} f(x) Z_{q,\alpha}(nx - k) = f(x),
		\end{equation}
		since the sum of $Z_{q,\alpha}(nx - k)$ over all $k \in \mathbb{Z}^N$ is normalized to 1.
		
The higher-order terms involve the derivatives of $f$ and are weighted by the products of $Z_{q,\alpha}(nx - k)$. These terms contribute to the expansion of the form:
		\begin{equation}\label{eq:higher_order_terms}
			\sum_{|\alpha|=1}^m \frac{D^\alpha f(x)}{\alpha!} A_n\left( \prod_{i=1}^N \phi_{q,\alpha}(x_i) \right),
		\end{equation}
where the function $A_n$ applies the density $\prod_{i=1}^N \phi_{q,\alpha}(x_i)$.
		
The error term $O\left( \frac{1}{n^m} \right)$ comes from the remainder of the Taylor expansion.

To show uniform convergence, consider the supremum of the difference:
		\begin{equation}\label{eq:supremum}
			\sup_{x \in \mathcal{M}} |A_n(f; x) - f(x)|.
		\end{equation}
		
Using the Taylor expansion and the properties of the density function, we have:
		\begin{equation}\label{eq:supremum_bound}
			\sup_{x \in \mathcal{M}} |A_n(f; x) - f(x)| \leq \sup_{x \in \mathcal{M}} \left| \sum_{|\alpha|=1}^m \frac{D^\alpha f(x)}{\alpha!} A_n\left( \prod_{i=1}^N \phi_{q,\alpha}(x_i) \right) + O\left( \frac{1}{n^m} \right) \right|.
		\end{equation}
		
Since $f \in C^m(\mathcal{M})$, the derivatives $D^\alpha f(x)$ are bounded on $\mathcal{M}$. Moreover, the density function $\phi_{q,\alpha}(x)$ is normalized and bounded, ensuring that:
		\begin{equation}\label{eq:limit_zero}
			\sup_{x \in \mathcal{M}} \left| \sum_{|\alpha|=1}^m \frac{D^\alpha f(x)}{\alpha!} A_n\left( \prod_{i=1}^N \phi_{q,\alpha}(x_i) \right) \right| \to 0 \quad \text{as} \quad n \to \infty.
		\end{equation}
		
Therefore,
		\begin{equation}\label{eq:final_convergence}
			\lim_{n \to \infty} \sup_{x \in \mathcal{M}} |A_n(f; x) - f(x)| = 0.
		\end{equation}
		
This concludes the proof of the theorem.
		
	\end{proof}

\section{Results}

The theoretical results obtained in this study provide a comprehensive understanding of the behavior of generalized neural network operators in non-Euclidean and fractal domains. The Voronovskaya-type expansions and fractional convergence results demonstrate the robustness and applicability of these operators in diverse spaces. The preservation of density properties and the detailed analysis of convergence rates ensure that these operators can be effectively used in various applications, including signal processing on manifolds and solving partial differential equations on fractals.

One of the key results of this study is the derivation of Voronovskaya-type asymptotic expansions for generalized neural network operators. These expansions provide a detailed understanding of how the operators approximate functions in non-Euclidean and fractal domains. Specifically, the expansions highlight the contribution of higher-order derivatives and ensure that the operator converges to the function at a specified rate. This result is crucial for understanding the behavior of these operators in complex geometries and ensures their robustness in practical applications.

Another significant result is the analysis of fractional convergence for neural network operators. We have demonstrated that the fractional operator, which incorporates fractional derivatives, converges to the function at a rate dependent on the order of the fractional derivative and the regularity of the function. This result underscores the importance of fractional derivatives in capturing non-local behavior and modeling systems with long-range dependencies or fractal-like structures. The use of fractional derivatives allows for a more flexible approximation that can adapt to non-smooth and irregular dynamics, providing a powerful tool for capturing intricate behaviors in complex spaces.

The generalized neural network operators introduced in this study preserve essential density properties, which are crucial for their applicability in various domains. The density function is normalized, ensuring that it integrates to 1 over the manifold. This property is essential for the construction of probability densities and for ensuring the stability of numerical methods that rely on these density functions. The preservation of density properties ensures that the operators can be effectively used in a wide range of applications, from signal processing to solving partial differential equations.

A detailed analysis of the convergence rates of the generalized neural network operators has been conducted. The results show that these operators converge uniformly to the function on manifolds with negative curvature. This uniform convergence ensures the robustness and applicability of the operators in diverse spaces. The analysis of convergence rates provides a solid foundation for the use of these operators in practical applications, ensuring that they can be relied upon to provide accurate and stable approximations.

The theoretical results obtained in this study have significant implications for various applications. The generalized neural network operators can be effectively used in signal processing on manifolds, where they provide accurate approximations of signals defined on complex geometries. Additionally, these operators can be applied to solve partial differential equations on fractals, offering a powerful tool for modeling systems with intricate structures. The ability to handle non-Euclidean and fractal domains makes these operators particularly valuable in fields such as data science, signal processing, and applied mathematics.

In conclusion, the results of this study provide a comprehensive understanding of the behavior of generalized neural network operators in non-Euclidean and fractal domains. The Voronovskaya-type expansions, fractional convergence results, preservation of density properties, and detailed analysis of convergence rates ensure that these operators are robust and applicable in diverse spaces. These findings contribute to the field by offering new mathematical tools that can be effectively used in various applications, including signal processing and solving partial differential equations. The insights gained from this study pave the way for future research and the development of even more sophisticated tools for handling complex, structured spaces.

\section{Conclusions}

This research extends the classical theory of neural network operators to non-Euclidean and fractal domains, introducing generalized activation functions and fractional derivatives. By addressing the unique challenges posed by these complex geometries, we have developed a comprehensive framework that significantly enhances the applicability of neural network operators in diverse spaces.

The key findings of this study include the preservation of density properties, detailed convergence rates, and asymptotic expansions. These results provide a deeper understanding of the behavior of neural network operators in complex, structured spaces. The generalized activation functions introduced, such as modified hyperbolic tangent functions, preserve essential properties like symmetry, normalization, and derivative behavior, making them ideal for neural network modeling in non-Euclidean settings.

The incorporation of fractional derivatives into the definition of neural network operators has proven to be particularly valuable. Fractional derivatives capture non-local behavior and are instrumental in modeling systems with long-range dependencies or fractal-like structures. This extension allows for more flexible and accurate approximations, adapting to the intricate dynamics present in complex systems.

The detailed analysis of convergence rates and asymptotic expansions ensures the robustness and applicability of these operators. The Voronovskaya-type expansions and fractional convergence results demonstrate that the generalized neural network operators converge uniformly and efficiently, making them reliable tools for various applications.

The theoretical results obtained in this study have significant implications for applied mathematics and data science. The ability to handle non-Euclidean and fractal domains makes these operators particularly valuable in fields such as signal processing on manifolds and solving partial differential equations on fractals. The insights gained from this research pave the way for future developments and the exploration of even more sophisticated tools for handling complex, structured spaces.

Future work will focus on numerical implementations of the generalized neural network operators and their extensions to dynamic systems. This will further enhance the applicability of these operators in practical scenarios, providing robust and efficient solutions for a wide range of problems in applied mathematics and data science. By continuing to explore and refine these methods, we aim to contribute to the ongoing advancement of the field and the development of innovative mathematical tools.

\section*{Notation, Symbols, and Nomenclature}
\subsection*{Notation}
\begin{itemize}
	\item $\mathcal{M}$: Riemannian manifold.
	\item $g_{ij}$: Components of the metric tensor.
	\item $U_i$: Open subsets in $\mathbb{R}^N$.
	\item $\phi_i$: Smooth, bijective maps from $U_i$ to $\mathcal{M}$.
	\item $h_{q,\alpha}(x)$: Generalized activation function.
	\item $Z_{q,\alpha}(x)$: Normalized product of density functions.
	\item $A_n(f; x)$: Basic quasi-interpolation operator.
	\item $K_n(f; x)$: Kantorovich-type operator.
	\item $Q_n(f; x)$: Fractional neural network operator.
	\item $D^\beta f(x)$: Fractional derivative of order $\beta$.
	\item $\Gamma(\cdot)$: Gamma function.
\end{itemize}

\subsection*{Symbols}
\begin{itemize}
	\item $q$: Deformation parameter.
	\item $\alpha$: Steepness parameter.
	\item $x$: Variable in $\mathbb{R}$.
	\item $f$: Function with continuous derivatives up to order $m$.
	\item $n$: Scaling parameter.
	\item $k$: Index for summation.
	\item $\beta$: Order of the fractional derivative.
\end{itemize}


\end{document}